\documentclass{article}
 \usepackage[ansinew]{inputenc}
 \usepackage[psamsfonts]{amssymb}
  \def\bk{{l\!k}}

\title{Graded Lie algebras with finite polydepth}
\author{Y. Felix, S. Halperin and J.-C.  Thomas}
\begin{document}
\maketitle

\begin{abstract}{ If $A$ is a graded connected algebra then we define a
new invariant, polydepth$\, A$, which is finite if
$\mbox{Ext}_A^*(M,A) \neq 0$ for some $A$-module $M$ of at most
polynomial growth. {\bf Theorem 1}: If $f : X \to Y$ is a continuous
map of finite category, and if the orbits of $H_*(\Omega Y)$ acting in
the homology of the homotopy fibre grow at most polynomially, then
$H_*(\Omega Y)$ has finite polydepth.
{\bf Theorem 5}:
If $L$ is a graded Lie algebra and polydepth $UL$ is finite then either $L$ is solvable and $UL$ grows at most
polynomially or else for some integer $d$ and all $r$,
$\sum_{i=k+1}^{k+d} \mbox{dim}\, L_i \geq k^r$, $k\geq$ some
$k(r)$.}
\end{abstract}

 \vspace{5mm}\noindent {\bf AMS Classification} : 55P35, 55P62,
 17B55

 \vspace{2mm}\noindent {\bf Key words} : Graded Lie algebras, loop
 space homology, depth,  solvable Lie algebras.

\vspace{5mm}

 We work over a field $\bk$ of characteristic different from 2. If $V = \{ V_k\}$ is a graded vector space we denote by $V^{\#} = \{\mbox{Hom}_{\bk}
 (V_k,\bk )\}$ the dual graded vector space. A graded Lie algebra is a
graded vector space $L$, equipped with a bilinear map $[\, ,\, ]$ $: L_i \times
L_j \to L_{i+j}$ satisfying
$$ [x,y] + (-1)^{ij} [y,x] = 0 $$
and
$$[x,[y,z]] = [[x,y],z] + (-1)^{ij}[y,[x,z]] $$
for $ x\in L_i$, $y\in L_j$, $z\in L$.
It follows that $3[x,[x,x]] = 0$ for $x$ of odd degree, and so if char $k=3$ we further require that $[x,[x,x]] = 0$. Finally we consider only graded Lie algebras satisfying $L=\{\, L_i\,\}_{i\geq 1}$ and each $L_i$ is finite dimensional. (Any graded vector space $V$ with each $V_i$ finite dimensional is said to have {\sl finite type}.)

The universal enveloping algebra of $L$ is denoted by $UL$ and it satisfies the classical Poincar\'e
Birkhoff Witt Theorem (in characteristic $3$ this uses the $[x,[x,x]]=0$ requirement).

Important examples  appear in   topology. Let $X$ be a  simply connected topological
space with rational homology of finite type. Then
the  rational homotopy Lie algebra $L_X$  of $X$ is defined by
$$ L_X = \pi_*(\Omega X)\otimes   {\mathbb Q}\,; \hspace{1cm} [\, ,\,] = \mbox{ Samelson
product}\,,$$
and the Hurewicz map extends to an isomorphism,  \cite{MM},
$$ UL_X \stackrel{\cong }{\longrightarrow} H_*(\Omega X; {\mathbb Q})\,.$$
Analogously, if $X$ is a finite $n$-dimensional and $r$-connected CW complex, then for primes
$p>n/r$, $H_*(\Omega X;F_p) = UE$ for some graded Lie algebra $E$ \cite{Hal}.

\vspace{3mm} If $M$ is a module over a (graded) algebra $A$ then the {\sl grade} of $M$, grade$\,
M$, is the least integer $q$ (or $\infty$) such that ${\rm Ext}_A^q(M,A)\neq 0$. And if
$V = \{\, V_i\,\}_{i\geq 0}$ is a graded vector space then $V$ has {\sl at most  polynomial growth}
if for some constant $C$, and some non-negative integer, $d$, $\sum_{i\leq n} {\rm dim}\, V_i \leq Cn^d$,
$n\geq 1$. In this case the least such $d$ is called the {\sl polynomial bound } for the growth of $V$
and is denoted by ${\rm polybd}\, (V)$. If $V$ does not have at most polynomial growth we put
${ \rm polybd}\, V = \infty$ and we say that $V$ grows faster than any polynomial.

In this paper we combine these two notions   in the

\vspace{2mm}\noindent {\bf Definition:} the {\sl polygrade} of an $A$-module, $M$, is the sum,
 ${\rm grade}\,M   + {\rm polybd}\, M$, and   the {\sl  polydepth of   $A$}
is the least integer (or $\infty$) occurring as the polygrade  of an $A$-module.

\vspace{2mm}
In the case $A=UL$ the unique augmentation  $UL \to \bk$ makes $\bk$  into a $UL$-module, and by
definition, the grade of $\bk$ is the {\sl depth} of $UL$. Since ${\rm polybd}\,\bk=0$ it follows that:
$$
 \mbox{polydepth}\, UL \leq \mbox{depth}\, UL
\eqno{(1)}
$$
Moreover (cf. Proposition 2) if ${\rm dim }\,L<\infty$ then equality holds. We shall abuse
notation  and refer to these invariants respectively as  polydepth$\, L$ and depth$L$.

\vspace{3mm}Note that ${\rm Ext}_{UL}^0(UL,UL)$ contains the identity map
and so
$$ {\rm polydepth}\, L \leq {\rm polybd}\,  UL \,. \eqno{(2)}$$

Observe as well that for any graded vector space $M$,
$\mbox{polybd}\, M = 0$ if and only if dim$\, M$ is finite. Thus
polydepth$\, L = 0$ if and only if depth$\, L = 0$, which happens
if and only if $L$ is finite dimensional and concentrated in odd
degrees.

Depth has been a useful concept in topology because, on the one hand, Lusternik-Schnirelmann
category satisfies \cite{modp}
$$ {\rm depth}\, L_X \leq {\rm cat}\, X$$
and, on the other hand (\cite{LAPG}, \cite{HAPG}, \cite{Engel}), finite depth has important
implications for the structure of a graded Lie algebra.

The purpose of this paper is to show that essentially the same implications follow from the
weaker hypothesis that polydepth$\, L$ is finite, while simultaneously identifying a   larger class
 of topological spaces and Lie algebras for which the weaker hypothesis holds.

Indeed, we have

\vspace{3mm}\noindent {\bf Theorem 1.}  {\sl If $F\to X \stackrel{f}{\longrightarrow}Y$ is a
fibration of path-connected spaces, then
$${\rm polydepth}\, H_*(\Omega Y) \leq   {\rm polybd} \,  H_*(F)  + {\rm cat}\, f\,.$$}

\vspace{2mm}\noindent {\bf proof:}  The fibration determines an action up to homotopy of
 $\Omega Y$ on $F$, which makes $H_*(F)$ into an $H_*(\Omega Y)$-module. According to
  \cite{grade}, grade$\, H_*(F) \leq $ cat$\, f$.
\hfill $\Box$

\vspace{3mm} Our main structural theorems read:

\vspace{3mm}\noindent {\bf Theorem 2.}  {\sl Let $E(L)$ denote the linear span of elements
 $x\in L_{\rm \scriptstyle even}$ such that ${\rm ad} x$ acts nilpotently on each $y \in L$.
 Then
$$ {\rm dim} \, E(L) \, \leq \, {\rm polydepth}\, L\,.$$}

\vspace{1.5cm}\noindent {\bf Theorem 3.} {\sl The following conditions on a graded Lie algebra $L$ are equivalent
\begin{enumerate}
\item[(i)] $L$ is the union of solvable ideals and ${\rm polydepth}\, L$ is
finite;
\item[(ii)] $UL$ grows at most polynomially (${\rm polybd}\, UL$ is
finite);
\item[(iii)] $L_{\rm \scriptstyle even}$ is finite dimensional, and for some constant $C$
$$ \sum_{i\leq n} {\rm dim}\, L_i \leq C \log_2 n\,, \hspace{1cm} n\geq 1\,.$$
\end{enumerate}
In this case $L$ is solvable.}

\vspace{3mm}\noindent {\bf Theorem 4.}  {\sl If $L$ is a graded
Lie algebra of finite polydepth then the union of the solvable
ideals of $L$ is a solvable ideal of finite polydepth.}

\vspace{3mm}\noindent {\bf Theorem 5.} {\sl Suppose ${\rm polydepth }\, L$ is finite and $L$ is not solvable. Then there is an integer $d$
such that for all $r\geq 1$:
$$ \sum_{i=k+1}^{k+d}{\rm dim}\, L_i \geq k^r\,, \hspace{1cm} k \geq \mbox{\rm some}
 \, k(r) \,.$$}

\vspace{3mm}\noindent {\bf Remark}. In \cite{growth} it is shown
that if $L = L_X$ where $X$ is a finite 1-connected CW complex,
then we may take $d=\mbox{dim}\, X$ in Theorem 5.

\section{Properties of polydepth}

\noindent {\bf Lemma 1:} {\sl If $M$ is a module for some graded algebra $A$ of finite type
and if ${\rm Ext}_A^q(M,A)\neq 0$ then ${\rm Ext}_A^q(A\cdot x,A) \neq 0$ for some $x$ in a
subquotient module of $M$.}

\vspace{3mm}\noindent {\bf proof:}  Recall that  $A^{\#} = Hom_{\bk}(A,\bk)$. Then ${\rm Ext}_A^q(M,A)$ is the
 dual of ${\rm Tor}_q^A(M,A^{\#})$ and a direct limit argument shows that for some
  $x_1, \ldots ,x_n \in M$, $${\rm Tor}_q^A(A\cdot x_1 + \ldots +A\cdot x_n,A^{\#}) \neq 0\,.$$
  Now use the exact sequence associated to the inclusion $A\cdot x_1 + \ldots + A\cdot x_{n-1}$
  in $A\cdot x_1 + \ldots + A\cdot x_n$.
\hfill $\Box$

\vspace{3mm}\noindent {\bf Corollary:}  Polydepth$\, A$ is the least $m$ such that
polygrade$\, N = m$ for some monogenic $A$-module $N$.

\vspace{3mm} \noindent {\bf Remark:} It follows from the Corollary that we may improve   Theorem
 1 to the inequality
$$ {\rm polydepth}\, H_*(\Omega Y) \leq {\rm polybd}\,
(H_*(\Omega Y)\cdot \alpha ) + {\rm cat}\, f\,, \hspace{5mm} {\rm some}\, \alpha \in H_*(F).
 \eqno{(3)} $$

 \vspace{3mm}\noindent {\bf Proposition 1.} {\sl Let $L$ be   a graded Lie algebra.
\begin{enumerate}\item[{\rm (i)}] Each ideal satisfies
   ${\rm polydepth}\, I \leq {\rm polydepth}\, L\,.$
\item[{\rm (ii)}] Let $E$ be a Lie subalgebra of   $L$. If $L$ has finite polydepth and
 if for each $x \in L/E$ the orbit  $UE\cdot x$ has at most   polynomial growth,
 then $E$ has finite polydepth.
\item[{\rm (iii)}] For $n$ sufficiently large the sub Lie algebra $E$ generated by
$L_{\leq n}$ satisfies ${\rm polydepth }\, E \leq {\rm polydepth}\, L$.
\end{enumerate}
}

\vspace{3mm}\noindent {\bf proof:}  (i) This follows from the Hochschild-Serre spectral
 sequence, converging from ${\rm Ext}_{UL/I}^p(\bk,{\rm Ext}_{UI}^q(M,UL))$ to ${\rm Ext}_{UL}^{p+q}(M,UL)$.
 (Note that since $UL$ is $UI$-free, grade$_{UI}(M)$ is the least $q$ such
 that ${\rm Ext}_{UI}^q(M,UL)\neq 0$.)

(ii) As in Lemma 1, ${\rm Ext}_{UL}^q(M,UL)$ is dual to ${\rm Tor}_q^{UL}(M,(UL)^{\#})$,
and this is the homology of the Cartan-Eilenberg-Chevalley complex $\land sL\otimes M \otimes (UL)^{\#}$.
 Write $L= E \oplus V$ and set $F_p =\land sE \otimes \land^{\leq p}sV\otimes M \otimes (UL)^{\#}$.
 This filtration determines a convergent spectral sequence, introduced by Koszul in \cite{Ko},
 and which is the Hochschild-Serre spectral sequence when $E$ is an ideal. The $E^1$-term of
 the spectral sequence is ${\rm Tor}_q^{UE}(\land^p sL/E \otimes M, (UL)^{\#})$, converging
 to ${\rm Tor}_{p+q}^{UL}(M,(UL)^{\#})$.

Each element $z\in \land^psL/E\otimes M$ is contained in a finite sum of $UE$-modules of the
 form $s(UE\cdot x_1) \land \cdots \land s(UE\cdot x_p) \otimes M$ and it follows that
  $${\rm polybd}\, (UE\cdot z) \leq p\cdot {\rm polybd}\, (UE\cdot x) + {\rm polybd}\, (M)$$
  for some $x\in L/E$. Choose $M$ so that   polydepth$\, L= $ polygrade$\, M$
  and apply Lemma 1 with $p+q = $ grade$\, M$.

(iii) If ${\rm Ext}_{UL}^p(M,UL)$ is non-zero and ${\rm polybd}\, (M) <\infty$ it suffices to choose $E$
 so that the restriction ${\rm Ext}_{UL}^p(M,UL) \to {\rm Ext}_{UE}^p(M,UL)$ is non zero.
\hfill $\Box$

\vspace{3mm}\noindent {\bf Corollary 1.} (of the proof of (ii)). {\sl  Suppose for some $k\geq 1$ that
 ${\rm polybd}\,(UE\cdot x)
 \leq k$, $x\in L/E$. Then ${\rm polydepth}\, E \leq k\, {\rm polydepth}\, L$.}

\vspace{3mm}\noindent {\bf Corollary 2.} {\sl Let $E$ be a sub-Lie algebra of a graded Lie algebra $L$.
 If $L$ has finite polydepth and $L/E$ has at most polynomial growth, then $E$ has finite polydepth.}

\vspace{3mm}\noindent {\bf Example 1.} Let ${\mathbb L}(V)$ be the free Lie algebra on a graded vector space
 $V$. Then for any graded Lie algebra $L$, $L\coprod {\mathbb
 L}(V)$ has depth 1. Thus
the injection $L \to
L\coprod L(V)$ shows that each graded Lie algebra is a sub-Lie
algebra of a Lie algebra of finite polydepth. The previous
corollary gives restriction on a Lie algebra $L$ for being a
sub-Lie algebra of a Lie algebra of finite polydepth, $K$, when
the quotient
has at most polynomial growth.

\vspace{3mm}\noindent {\bf Proposition 2.} {\sl If $L$ is a finite dimensional graded Lie algebra then
$${\rm polydepth}\, L = {\rm depth}\, L\,.$$
}

  \noindent {\bf proof:} As observed in the introduction, ${\rm polydepth}\,
 L \leq {\rm depth}\, L$. On the other hand, by Lemma 1, ${\rm polydepth}\, L = {\rm polygrade}\,
 M$ for some monogenic module $M = UL\cdot x$. Now Theorem 3.1 in \cite{LAPG} asserts
 that ${\rm polygrade}\, M = {\rm depth}\, L$.
\hfill $\Box$

\section{Proof of Theorem 2}

Suppose $I\subset L$ is an ideal.
If ${\rm Ext}_{UL}^m(M,UL) \neq 0$, then  ${\rm Ext}_{UL/I}^p({\rm Tor}_q^{UI}(M,k),UL/I)\neq
0\,,$
some $p+q =m$. (Same proof as in: \cite{LAPG}, Lemma 4.3, for the case $M=\bk$).
By Lemma 1 there is a monogenic $UL/I$-module $N$ such that $N$ is a subquotient of $\mbox{Tor}_q^{UI}(M,\bk)$,
and
grade $N\leq p$.

Now suppose   $L/I$ is finite dimensional. Then   Theorem 3.1 in \cite{LAPG} asserts that
$$ {\rm grade}\,N + {\rm polybd}\, N = {\rm dim}\, (L/I)_{\rm \scriptstyle even}\,.$$
On the other hand, write $(L/I)_{\rm \scriptstyle even} = V \oplus W$ where $V$ is
the image of $E(L)$.
Let $x_i\in L_{\scriptsize even}$, $y_j\in L_{\scriptsize odd}$ and $z_k\in E(L)$ represent
 respectively bases of $W$, $(L/I)_{\scriptsize odd}$ and $V$. Then the elements
 $x_1^{k_1}\cdots x_s^{k_s}y_1^{\varepsilon_1}
\cdots y_t^{\varepsilon_t}z_1^{m_1}\cdots z_u^{m_u}$, where  $\varepsilon_i = 0$
or $1$, represent a basis for $UL/I$.  Choose the $z_k$ to act locally nilpotently in $L$.
Then this basis applied to any
 $\omega \in \land^qsI$, shows that
${\rm polybd}\,(UL/I)\cdot \omega \leq {\rm dim} W$.
Hence if $u \in \land^qsL\otimes M$ represents a generator of $N$ then
$${\rm polybd}\, N \leq {\rm polybd}\, (UL/I\cdot u)
\leq {\rm polybd}\, M+ {\rm dim} W\,.$$
Substitution in the equation above gives
$$
\renewcommand{\arraystretch}{1.5}
\begin{array}{ll}
{\rm dim}\,(L/I)_{\rm \scriptstyle even} & \leq {\rm grade}\, N +{\rm polybd}\, M+ {\rm dim}
\, W \\
& \leq {\rm grade}\, M + {\rm polybd}\, M + {\rm dim} \, W\,.
\end{array}
\renewcommand{\arraystretch}{1}
$$
Choose $M$ so grade$\, M + {\rm polybd}\, (M) = {\rm polydepth}\, L$ and choose $I = L_{> 2k}$.
Then $V \cong E(L)_{\leq 2k}$ and we have
$${\rm dim}\, E(L)_{\leq 2k} \leq {\rm polydepth}\, L\,.$$
Since this holds for all $k$ the theorem is proved.
\hfill $\Box$

\section{Solvable Lie algebras}

 \vspace{3mm}\noindent {\bf Lemma 2.} {\sl Let $L$ be a Lie algebra concentrated in odd
 degrees.
 Then ${\rm Ext}_{UL}(-,UL) = {\rm Hom}_{UL}(-,UL)$. In particular
$${\rm polydepth }\, L = {\rm polybd}\, UL\,.$$}

 \noindent {\bf proof:}  Since $L = L_{\rm \scriptstyle odd}$   it is necessarily abelian.
 Now ${\rm Ext}_{UL}(-,UL)$ is the dual  of ${\rm Tor}^{UL}(-,(UL)^{\#})$ and this is the
 limit of ${\rm Tor}^{UL_{\leq n}}(-,(UL)^{\#})$, which dualizes
 to ${\rm Ext}_{UL_{\leq n}}(-,UL)$. Since $UL_{\leq n}$ is a finite dimensional exterior algebra
 and $UL$ is $UL_{\leq n}$-free it follows that ${\rm Ext}^+_{UL_{\leq n}}(-,UL) =
 0$, and so $\mbox{Ext}^+_{UL}(-, UL) = 0$.

Finally, since ${\rm Ext}^0_{UL}(UL,UL)$ is non-zero, ${\rm polydepth}\, L \leq {\rm polybd}\,
UL$.
On the other hand if ${\rm polydepth}\, L = m<\infty$, then for some $M$,
we have ${\rm Ext}_{UL}^p(M,UL)\neq 0$ and ${\rm polybd}\, M = m-p$. By the above,
$p=0$ and so there is a non zero $UL$-linear map $f : M \to UL$. Any $f(m)$ is in some
 $UL_{<n}$ and if $f(m)\neq 0$ it follows that
 $UL_{\geq n} \stackrel{\cong}{\longrightarrow} UL_{\geq n}\cdot m$. This implies ${\rm polybd}\, M
 \geq {\rm polybd}\, UL$ and ${\rm polydepth}\, L \geq {\rm polybd}\, UL$.
\hfill $\Box$

\vspace{3mm}\noindent {\bf Lemma 3.}  {\sl Let $L$ be a graded Lie algebra of finite polydepth. If
  $I$ is an ideal in $L$ and
${\rm polybd} \,  I  <\infty$ then ${\rm polydepth}\, L/I < \infty$.}

\vspace{3mm}\noindent {\bf proof:}  Choose $M$ so that ${\rm polygrade}\, M = {\rm polydepth}\, L$.
If $m = {\rm grade}\, M$ then it follows (as in \cite{LAPG}, proof of Theorem 4.1 for the case $M=\bk$)
that for some $p$,
$${ \rm Ext}_{UL/I}^p({\rm Tor}_{m-p}^{UI}(M,\bk),UL/I) \neq 0\,.$$
Since ${\rm Tor}_{m-p}^{UI}(M,\bk)$ is a subquotient of $\land^{m-p}sI\otimes M$ it follows that it has
polynomial growth at most equal to $(m-p) {\rm polybd}\, I $
\hfill $\Box$

\vspace{3mm}\noindent {\bf proof of Theorem 3:}

(i) $\Rightarrow$ (ii).  Let $I$ be the sum of the ideals in $L$ concentrated in odd degrees. Then
 $I$ is an ideal of this form, necessarily abelian, and $L/I$ has no ideals concentrated in odd degrees.
 Moreover ${\rm polybd}\, UI ={\rm polydepth}\, I \leq {\rm polydepth}\, L$ (Lemma 2 and Proposition 1)
 and hence ${\rm polydepth}\, L/I < \infty$ (Lemma 3).

Next we show that every solvable ideal $J$ in $L/I$ is finite dimensional, by induction on the solvlength.
Indeed, if $J$ is abelian then $J_{\rm \scriptstyle even} = E(J)$. Since
${\rm polydepth}\, J \leq {\rm polydepth}\, L/I$ (Proposition 1), Theorem 2 asserts
that $J_{\mbox{\scriptsize even}}$
 is finite dimensional. Thus for some $r$, $J_{\geq r}$
is an ideal concentrated in odd degrees; i.e. $J_{\geq r} = 0$.

Now if $J$ has solvlength $k$ then its $(k+1)$st derived algebra is abelian and so finite dimensional.
Thus for some $r$, $J_{\geq r}$ has solvlength $k-1$. By induction, $J$ is finite dimensional.

By hypothesis $L/I$ is the sum of its solvable ideals. Since these are finite dimensional, each
$x\in (L/I)_{\rm \scriptstyle even}$ acts locally nilpotently. Thus
$(L/I)_{\rm \scriptstyle even} = E(L/I)$, and this is finite dimensional by Theorem 2.
But $L_{\rm \scriptstyle even} \cong (L/I)_{\rm \scriptstyle even}$ since $I$ is concentrated in odd degrees.

Suppose $L_{\rm \scriptstyle even} \subset L_{\leq 2n}$. Since $L_{>2n}$ is an ideal in
odd degrees of finite polydepth, ${\rm polybd}\, UL_{>n} <\infty$,  while trivially
${\rm polybd}\, UL/L_{>2n}<\infty$. Hence ${\rm polybd}\, UL < \infty$.

\vspace{2mm} (ii) $\Rightarrow$ (iii).  Clearly ${\rm polybd}\, UL \geq {\rm dim}
L_{\rm \scriptstyle even}$, so the latter must be finite. It is trivial from the Poincar\'e
 Birkhoff Witt theorem that if $\sum_{i\leq n}{\rm dim}\, L_i = d(n)$ then
 $$\sum_{i\leq nd(n)} {\rm dim}\, (UL)_i \geq 2^{d(n)}\,.$$ Thus $2^{d(n)} \leq K{[nd(n)]^r} $
 for some constant $K$ and some integer $r$, $r\geq 1$.  It follows that
 $d(n) \leq \log_2K + r\log_2n+r\log_2d(n) \leq   r\log_2n + \frac{1}{2} d(n)$, $n$
 sufficiently large.

\vspace{2mm} (iii) $ \Rightarrow$ (i).  Choose $N$ so that $I = L_{\geq N}$ is
concentrated in odd degrees. Then $UI$ is an exterior algebra and so
$$ \sum_{i\leq n} {\rm dim}\, (UL)_i \leq 2^{ \sum_{i\leq n} {\rm dim}\, L_i} \leq n^C$$
for some constant $C$.
Thus,
since $L/I$ is finite dimensional, ${\rm polybd}\, UL$ is finite. The identity of $UL$ is
in ${\rm Ext}^0_{UL}(UL,UL)$ and so ${\rm polydepth}\, L \leq {\rm polybd}\, UL$.

Finally, since $I$ is abelian and $L/I$ is finite dimensional, $L$ itself is solvable. This
also proves the last assertion.
\hfill $\Box$

\vspace{3mm}\noindent {\bf proof of Theorem 4:}  This is immediate
from Proposition 1(i) and Theorem 3. \hfill $\Box$

\vspace{3mm}\noindent {\bf Proposition 3.} {\sl Suppose $I$ is a solvable ideal in a Lie algebra
$L$ of finite polydepth. Then
\begin{enumerate}
\item[(i)] ${\rm polydepth}\, L/I \leq {\rm polydepth}\, L$
\item[(ii)] ${\rm dim} I_{\rm \scriptstyle even} \leq {\rm polydepth}\, I \leq {\rm polybd}\, UI$.
\end{enumerate}}

\vspace{3mm}\noindent {\bf proof:}  (i)   As noted in the proof of Lemma 3,
$${\rm Ext}^p_{UL/I}({\rm Tor}_{m-p}^{UI}(M,\bk),UL/I)\neq 0\,,$$
where $m + {\rm polybd}\, M = {\rm polydepth} \, L$. Also the Tor is a subquotient
of
$\land^{m-p}sI\otimes M$. By Theorem 3
$$ \sum_{i\leq n} {\rm dim}\, (\land^{m-p}sI\otimes M)_i \leq  \left( C_1\log_2 n \right)^{m-p}\cdot
C_2\cdot
n^{{\rm polybd}\,
M}\,.$$
Hence $\mbox{polybd}\, \mbox{Tor}_{m-p}^{UI}(M,\bk) \leq
\mbox{polybd}\, M + 1$, and so $\mbox{polydepth}\, L/I \leq
\mbox{polydepth}\, M + p + 1$.

If $p<m$ then this gives $\mbox{polydepth}\, L/I \leq
{\mbox{polydepth}}\, L$. If $p=m$ then
$\mbox{Tor}_{m-p}^{UI}(M,\bk) = M \otimes_{UI}\bk$. Hence in this
case $\mbox{polybd}\, \mbox{Tor}_{m-p}^{UI}(M,\bk ) \leq
\mbox{polybd}\, M$ and again $\mbox{polydepth}\, L/I \leq
\mbox{polydepth}\, L$.

\vspace{1mm}  (ii) Since $I_{\rm \scriptstyle even} \subset I_{\leq 2n}$, some $n$ (Theorem 3)
we may apply the first assertion to obtain
$$
\begin{array}{ll}
{\rm dim} \, I_{\rm \scriptstyle even}   & =  {\rm depth}\, I/I_{>2n}
\,\,\,{\rm (cf. \cite{FHJLT})}\\
&={\rm polydepth}\, I/I_{>2n}\,\,\, {\rm (Proposition \, 2)}\\
& \leq {\rm polydepth} I\,.
\end{array}
$$
The second inequality has already been observed: ${\rm polybd} \, UI = {\rm polygrade}\,
UI \geq {\rm polydepth}\, I$.

\hfill $\Box$

  \vspace{3mm} \noindent{\bf Example 3.} Consider a Lie algebra $L$ concentrated in odd
degrees with a basis $\{\,x_i\,\, i\geq 1\,\}$ satisfying the degree
relations
$$
\mbox{deg} x_i > \sum_{j<i} \mbox{deg} x_j \,.$$
Then for each $n$, ${\rm dim }(U L)_n \leq 1$. The identity on $UL$ shows
that $\mbox{polydepth}\, L = 1$.

\vspace{3mm}\noindent {\bf Example 4.}   Consider the graded Lie algebra $L = {\mathbb L}(a,x_n)_{n\geq 2}/I$,
with $\mbox{deg} a = 2$, $
\mbox{deg} x_n = 2^n+1$, and where $I$ is generated by the relations
$$ [({\rm ad} a)^k x_r,({\rm ad} a)^lx_s]=0, \hspace{3mm} k,l\geq 0, r,s\geq 2\,, \hspace{5mm}{\rm and} \hspace{5mm} {\rm ad}^{n+1}(a)(x_n) = 0\,.$$
Then   $\mbox{polybd}\, UL = 2$,   so that $L$ has finite polygrade.
On the other hand, $L$ is solvable but not
nilpotent, and is the union of the infinite sequence of the finite dimensional Lie algebras
$I_N$ generated by $a, x_2, \ldots, x_N$.

\vspace{3mm}\noindent {\bf Proposition 4.} {\sl Let $L $ be
the direct sum of non-solvable Lie algebras $L(i)$, $i\leq n$. If polydepth$\, L(i) < \infty$  for
$1\leq i\leq n$, then
$$n \leq {\rm polydepth}\, L \leq \sum_i{\rm polydepth }\,
L(i)\,.$$}

\vspace{3mm}\noindent {\bf proof.} We first prove by induction on $n$ that for any $UL$-module $M$
that has at most polynomial growth, we have ${\rm Ext}_{UL}^{<n}(M,UL) = 0$. Consider the
 Hochschild-Serre spectral sequence
$$ {\rm Ext}_{UL(1)}^p({\rm Tor}_q^{U(L(2)\oplus \ldots \oplus L(n))}(M,(U(L(2) \oplus \ldots
L(n))^{\#}), UL(1)) \Rightarrow {\rm Ext}_{UL}^{p+q}(M,UL) \,.\eqno{(4)}$$
Since $L(1)$ commutes with the other $L(i)$ it follows that for each monogenic $UL(1)$-module $N$ that is a
subquotient
of $\mbox{Tor}_q^{U(L(2) \oplus \cdots \oplus L(n))} (M, U(L(2)\oplus \cdots \oplus L(n))^{\#})$
we have $\mbox{polybd}\, N \leq \mbox{polybd}\, M$.

Now since $L(1)$ is not solvable, $\mbox{polybd}\,UL(1)=\infty$   and the argument in the proof of Lemma
2 shows that $\mbox{Ext}^0_{UL(1)}
(N,UL(1)) = 0$. Thus (Lemma 1) the left hand in (4) vanishes for $p=0$. By induction on $n$ it vanishes for $q<n-1$ and
so $\mbox{Ext}_{UL}^{<n}(M,UL) = 0$. Thus $\mbox{polydepth}\, L\geq n$.

On the other hand,   there are
$UL(i)$-modules $M(i)$ such that $\mbox{polygrade}\, M(i) =
\mbox{polydepth}\, L(i)$. Then $\otimes_{i=1}^n M(i)$ is a
$UL$-module that has at most polynomial growth and whose polygrade
is the sum of the polygrades of the $M(i)$.
\hfill $\Box$

\section{Growth of Lie algebras}

\vspace{3mm}\noindent {\bf Proposition 5.}  {\sl Let $L$ be a non
solvable graded Lie algebra of finite polydepth. Then for each
integer $r\geq 1$ there is a positive integer $d(r)$ such that
$$\sum_{i=k+1}^{k+d(r)} \mbox{dim}\, L_i \geq k^r\,, \hspace{1cm}
k \, \mbox{sufficiently large}\,.$$}

\vspace{3mm}\noindent {\bf proof:}  We distinguish two cases.

\noindent {\bf Case
A: } $L_{\mbox{\scriptsize even}}$ contains an infinite dimensional abelian sub Lie algebra
$E$.

Choose $n$ so that dim $E_{\leq n} \geq (r+3)
\mbox{polydepth}\, L$.  Then there is a finite sequence
$$L = I(0) \supset I(1) \supset \cdots \supset I(l)$$
in which $I(j)$ is an ideal in $I(j-1)$ and $I(l)_{\leq n} =
E_{\leq n}$.

By Proposition 1, polydepth$\, I(q) \leq$ polydepth$\, L$. Thus
without loss of generality we may suppose that $L = I(l)$; i.e.
that $L_{\leq n}$ is an abelian sub Lie algebra concentrated in
even degrees and that dim $L_{\leq n} \geq (r+3)\mbox{polydepth}\,
L$.

Let $M$ be a $UL$-module such that $\mbox{grade}\, M +
\mbox{polybd}\, M = \mbox{polydepth}\,L$ and put $m =
\mbox{grade}\, M$. As observed in the proof of Proposition 1(ii),
$\mbox{Ext}_{UL_{\leq n}}^q(\land^psL_{>n}\otimes M, UL_{\leq
n})\neq 0$, for some $p+q = m$. It follows that for some $z \in
\land^psL_{>n}\otimes M$,
$$\mbox{polybd}\, UL_{\leq n}\cdot z + q \geq \mbox{dim}\, L_{\leq
n}$$
(Theorem 3.1 in \cite{LAPG}). Hence for some $x\in L$,
$$p(\mbox{polybd} \, UL_{\leq n}\cdot x) \geq \mbox{dim}\, L_{\leq
n} - q - \mbox{polybd}\, M\,.$$

Since $p+q+ \mbox{polybd} \,M = \mbox{polydepth}\, L$ we conclude
that $$  (2 + \mbox{polybd}\, UL_{\leq n}\cdot
x) \cdot \mbox{polydepth}\, L \geq \mbox{dim}\, L_{\leq n}\geq (r+3)\mbox{polydepth}\, L\,.$$
As observed in the introduction, since dim$\, L = \infty$,  polydepth$\, L>0$. It follows that
$$\mbox{polybd}\, \left( UL_{\leq n}\cdot x\right)\geq r+1\,. \eqno{(5)}$$

On the other hand, $UL_{\leq n}$ is the polynomial algebra $\bk
[y_1, \ldots , y_s]$ on a basis $y_1, \ldots , y_s$ of $L_{\leq
n}$. Because of (5) it is easy to see (induction on $s$) that this
basis can be chosen so that for some $w\in UL_{\leq n}\cdot x$,
$\bk[y_1, \ldots , y_{r+1}] \to \bk [y_1, \ldots , y_{r+1}]\cdot
w$ is injective. Put $d = \prod_{i=1}^{r+1} \mbox{deg}\, y_i$ and
note that
$$\sum_{i=k+ \mbox{\scriptsize deg} w + 1}^{k+ \mbox{\scriptsize deg}\, w + d}
\mbox{dim}\, L_i \geq \sum_{i=k+1}^{k+d} \mbox{dim}\,\bk [y_1, \ldots ,
y_{r+1}]_i \geq \frac{1}{r!} k^r\,.$$
It follows that for any $r$, and $k$ sufficiently large,
$$\sum_{i=k+1}^{k+d} \mbox{dim}\, L_i \geq k^{r-1}\,.$$
This proves Proposition 5 in case A.

\vspace{3mm}\noindent {\bf Case B}  Every abelian sub Lie algebra
of $L_{even}$ is finite dimensional.

Let $I$ be the sum of the solvable ideals in $L$. Then
$I_{\mbox{\scriptsize even}}$ is finite dimensional and $\mbox{polydepth}\,
L/I$ is  finite (Theorem 3 and Proposition 3). Thus all abelian
sub Lie algebras of $L/I$ are finite dimensional. Thus it is
sufficient to prove Case B when $L$ has no solvable ideals.

There are now two possibilities: either $L=L_{\mbox{\scriptsize
even}}$, or $L$ has elements of odd degree. In the latter case the
sub Lie algebra generated by $L_{\mbox{\scriptsize odd}}$ is an
ideal, hence non-solvable and of finite polydepth. Let $L(s)$
denote the sub Lie algebra generated by the first $s$ linearly
independent elements $x_1, \ldots , x_s$ of odd degree. For $s$
sufficiently large, $\mbox{polydepth}\, L(s) \leq
\mbox{polydepth}\, L$ (Proposition 1) and $\mbox{dim}\,
L(s)_{\mbox{\scriptsize even}} > \mbox{polydepth}\, L$ (obvious).
Thus $L(s)$ cannot be solvable (Theorem 3). In other words, we may
assume that either $L=L_{\mbox{\scriptsize even}}$ or else $L$ is
generated by finitely many elements $x_1, \ldots , x_s$ of odd
degree. In either case set $E = L_{\mbox{\scriptsize even}}$, and note that $\mbox{dim}\, E$ is infinite.

Define a sequence of
elements $z_i$ and sub Lie algebras $E(i)$ by setting $E(1) = E$,
$z_i $ is a non-zero element in $E(i)$ and $E(i+1) \subset E(i)$
is the sub Lie algebra of elements on which ad$\, z_i$ acts
nilpotently.

Since $E$ contains no infinite dimensional abelian Lie algebra
some $E(N+1) = 0$ and $E(1)/E(2) \oplus \cdots \oplus E(N)/E(N+1)$
is a graded vector space isomorphic with $E$.

Put $d = \prod \mbox{deg}\, z_i$ and $d_i = d/\mbox{deg}\, z_i$.
Then $(\mbox{ad}\, z_1)^{d_1} \oplus \cdots \oplus
(\mbox{ad}\,z_N)^{d_N}$ is an injective transformation of
$E(1)/E(2) \oplus \cdots \oplus E(N)/E(N+1)$ of degree $d$. Since
this space is isomorphic with $E$ it follows that
$$\sum_{i=k+1}^{k+d}\mbox{dim}\, E_i \geq \frac{d}{k+d}
\sum_{i=1}^{k+d} \mbox{dim}\, E_i\,, \hspace{1cm} k\geq 1\,.\eqno{(6)}$$

On the other hand, choose $n$ so that
$$\mbox{dim}\, E_{\leq n} \geq (r+3)\cdot\mbox{polydepth}\, L\,.$$
(This is possible because $E$ is infinite dimensional.) Set $I =
L_{>n}$.
Let $M$ be
an $L$-module with $\mbox{polydepth}\, L = \mbox{polygrade}\, M$.
As in the proof of Theorem 4.1 in \cite{LAPG},
$\mbox{Ext}^p_{UL/I} (\mbox{Tor}_q^{UI}(M,\bk),UL/I)\neq 0$ for
$p+q = \mbox{grade}\, M$. Thus $p,q$ and $\mbox{polybd}\, M$ are
all bounded above by $\mbox{polydepth}\, L$.
Now Theorem
  3.1 of \cite{LAPG} asserts that for some $\alpha
\in \mbox{Tor}_q^{UI}(M,\bk)$, $UL/I\cdot \alpha$ has polynomial
growth at least equal to $\mbox{dim}(L/I)_{\mbox{\scriptsize even}}
- p $.  This means that for some
positive $C$,
$$\sum_{i\leq k} \mbox{dim}\,(UL/I\cdot \alpha)_i \geq
Ck^{(\mbox{\scriptsize dim}\, (L/I)_{\mbox{\scriptsize even} })-p}\,,\hspace{1cm}
k\,\,
\mbox{sufficiently   large}\,.$$
Since $\mbox{Tor}_*^{UI}(M,\bk)$ is the homology of
$\land^*sI\otimes M$, it follows that
$$\sum_{i\leq k} \mbox{dim}\, \mbox{Tor}_q^{UI}(M,\bk)_i \leq
(\sum_{i\leq k}\mbox{dim}\, L_i)^q \sum_{i\leq k}\mbox{dim}\,
M_i\,.$$
But $(L/I)_{\mbox{\scriptsize even}} \cong E_{\leq n}$ and so a
quick calculation gives
$$\sum_{i\leq k}\mbox{dim}\, L_i \geq Kk^{r+1}\,, \hspace{1cm}
\mbox{$k$ sufficiently large}\,. \eqno{(7)}$$

Finally, recall that either $L=L_{\mbox{\scriptsize even}}$ or
else $L$ is generated by the elements of odd degree $x_i$. In the
former case $L = E$ and the Proposition follows from (6) and (7).
In the second case we have $L_{\mbox{\scriptsize odd}} = [x_1,E] +
\cdots + [x_s,E]+ \bk x_1 + \cdots + \bk x_s$, and hence (7)
yields
$$\sum_{i\leq k} E_i \geq \frac{K}{s+1} k^{r+1} + s\,,\hspace{1cm}
\mbox{$k$ sufficiently large}\,.$$
Combined with (6) this formula gives the Proposition. \hfill
$\Box$

\vspace{3mm}\noindent {\bf Proof of Theorem 5:}  Since $L$ is not
solvable we may choose $n$ so that
$\mbox{dim}\,(L_{\mbox{\scriptsize even}})_{\leq n} >
\mbox{polydepth}\, L$ (Theorem 3), and so that the sub Lie algebra
generated by $L_{\leq n}$ satisfies $\mbox{polydepth}\, E\leq
\mbox{polydepth}\, L$ (Proposition 1).  Then $E$ is not solvable
(Proposition 3).

Let $x_1, \ldots , x_s$ generate $E$, and put $d = \mbox{max
deg}\, x_i$. Letting $UE$ act via the adjoint representation on
$E$ we have that
$$UE_{[0,q]}\cdot E_{[k+1,k+d]} \supset E_{[k+1,k+q]}\,.$$
For any $r\geq 1$ choose $q = q(r)$ so that $\sum_{i=k+1}^{k+q}
\mbox{dim}\, E_i \geq k^{r+1}$, $k$ sufficiently large
(Proposition 5). Then
$$\sum_{i=k+1}^{k+d} \mbox{dim}\, L_i \geq \sum_{i=k+1}^{k+d} \mbox{dim}\,
E_i\geq \frac{1}{\mbox{dim}\, UE_{[0,q]}}k^{r+1} \geq k^r\,,
\hspace{1cm} \mbox{$k$ sufficiently large}\,.$$
Since $d$ is independent of $r$, the Theorem is proved. \hfill
$\Box$

\vspace{5mm} Institut de Math\'ematiques,

Universit\'e Catholique de Louvain

1348 Louvain-La-Neuve, Belgium

\vspace{5mm} College of Computer, Mathematical and Physical Sciences

University of Maryland

College Park, MD 20742-3281, USA

\vspace{5mm} Facult\'e des Sciences

Universit\'e d'Angers

Bd. Lavoisier, 49045 Angers, France

\end{document}